\RequirePackage{fix-cm}
\documentclass[smallextended]{svjour3}       
\smartqed  
\usepackage{graphicx}
\usepackage{amssymb}
\usepackage{amsmath}
\usepackage{makeidx}

\newtheorem{thm}{Theorem}

\newtheorem{cor}{Corollary}

\newtheorem{prop}{Proposition}
\newcommand{\finish}{\hfill$\Box$\vspace{0.2cm}}

\newcommand{\prf}{\noindent{\bf Proof:\ }}

\newcommand{\E}{{\rm I \!E}}
\newcommand{\p}{{\rm I \!P}}

\begin{document}

\title{Maximum Loss and Maximum Gain of Spectrally Negative L\'{e}vy Processes \thanks{This work is supported by the Scientific and Technological Research Council of Turkey, TUBITAK Project No.110T674}}



\author{Ceren Vardar-Acar \and Mine \c{C}a\u{g}lar}


\institute{Ceren Vardar Acar \at
              Department of Statistics, Middle East Technical University, Ankara, Turkey \\
              \email{cvardar@metu.edu.tr}   \\
           \and
            Mine \c{C}a\u{g}lar\at
            Department of Mathematics, Ko\c{c} University, Istanbul, Turkey \\
              \email{mcaglar@ku.edu.tr}
           }

\date{Received: date / Accepted: date}
\maketitle

\begin{abstract}
The joint distribution of the maximum loss and the maximum gain is obtained for a spectrally negative L\'evy process until the passage time of a given level. Their marginal distributions up to an independent exponential time are also provided. The existing formulas for Brownian motion with drift are recovered using the particular scale functions.
\end{abstract}

\section{Introduction}
\label{intro}

The maximum loss, or maximum drawdown of a process $X$ is the supremum of $X$ reflected at its running supremum.  The  motivation comes from mathematical finance as it is useful to quantify the risk associated with the performance of a stock. The loss process has been studied for Brownian motion (Salminen and Vallois, 2007, Vardar-Acar et al. 2013), and some L\'evy processes (Mijatovic and Pistorius, 2012).
The maximum loss at time $t>0$ is formally defined by
\[
M_t^ - : = \mathop {\sup }\limits_{0 \le u \le v \le t} (X_u^{} - X_v^{}) \; ,
\]
which is equivalent to
 $\mathop {\sup }\limits_{0 \le v \le t} (\mathop {\sup }\limits_{0 \le u \le v} ({X_u} - {X_v}))$ and $ \mathop {\sup }\limits_{0 \le v \le t} ({S_v} - {X_v})$, that is, the supremum of the reflected process $S-X$, or the so-called loss process, where $S$ denotes the running supremum.
Similarly, a counterpart quantity called the maximum gain is given by $M_{t}^ +  := \mathop {\sup }\limits_{0 \le u \le v \le t} (X_v^{} - X_u^{})= \mathop {\sup }\limits_{0 \le v \le t} ({X_v} - {I_v})$ where $I$ is the running infimum. It is easy to see that the maximum gain of $X$ is the maximum loss of the process $-X$.

A spectrally negative L\'evy process $X$ is a L\'evy process with no positive jumps, that is,   its L\'evy measure $\Pi$ is concentrated on $(-\infty, 0)$. Let $\psi$ denote the Laplace exponent of $X$, as given by
$
\E(\exp \lambda {X_t})  = \exp ( t\,\psi (\lambda ))
$,
which is finite for all $\lambda \geq 0$ (Bertoin, 2007, pg.188). It is related to the characteristic exponent $\Psi$ by $\psi(\lambda)  = -\Psi (-i\lambda )$.  The Laplace exponent of a spectrally negative L\'evy process is given by
\[
\psi (\lambda ) =  - \mu \lambda  + \frac{{{\sigma ^2}}}{2}{\lambda ^2} + \int\limits_{( - \infty ,0)}^{} {({e^{\lambda x}} - 1 - \lambda x\;1_{\{x > -1\}})\, \Pi (dx)}
\]
where $\mu \in \mathbb{R}$, and $\sigma\neq0$. Note that Brownian motion with drift $\mu$ can be taken as a special case when $\Pi\equiv 0$.

In this paper, we first consider the joint distribution of the maximum loss and the maximum gain of a spectrally negative L\'evy process up to the passage time above a level $\beta$, denoted by $\tau_\beta$. It is found as
\[\p\{M_{{\tau_\beta }}^ -  \leq u,M_{{\tau_\beta }}^ +  \leq v\} = \left\{ \begin{array}{ll}
 e^{ - \beta W'_+(u)/W(u)}{\rm{                            }}& u \le v - \beta  \\
 \frac{{W(v - \beta )}}{{W(u)}}e^{ - (v - u){\kern 1pt} W'_+(u)/W(u)}{\rm{        }}& v - \beta  \le u \le v \\
 \frac{{W(v - \beta )}}{{W(v)}}{\rm{                         }} & v \le u\, \\
 \end{array} \right.\]
for $u\ge 0$ and $v\ge \beta$, where $W$ is the so-called scale function, which is defined as the inverse Laplace transform of $1/\psi(\lambda)$. We recover the analogous result as given in Salminen and Valois (2007) for Brownian motion with drift $\mu$ and $\sigma=1$, using in particular $W(x)=(1-e^{-2\mu x})/{2\mu}$.

Second, we find the marginal distributions of the maximum loss and the maximum gain up to an independent exponential time by drawing upon the existing results by Avram et al. (2004) and Pistorius (2004). Let $\gamma>0$ denote the parameter of the exponential random variable. The formulas obtained are
 \begin{eqnarray*}
\p\{M_T^- >a\}& = & {Z^{(\gamma)}}(a) - \gamma\, {[{W^{(\gamma)}}(a)]^2}/W_ + ^{(\gamma)'}(a) \\
\p\{M_T^+ >a\}& = & \frac{{Z^{(\gamma)}}(0)}{{Z^{(\gamma)}}(a)}
 \end{eqnarray*}
for $a>0$, where $W^{(\gamma)}$ and $Z^{(\gamma)}$ are the so-called $\gamma$-scale functions (Kyprianou, 2014). We show that these expressions yield the available results for the special case of Brownian motion with drift when the corresponding $\gamma$-scale functions are used.

The paper is organized as follows. In  Section 2, we give the essential definitions related to  the excursions and the scale function of a spectrally negative L\'evy process. In Section 3, the joint distribution of the maximum  loss and the maximum gain are obtained up to a time the process passes above a given level.  The marginal distributions  are found up to an independent exponential time in Section 4.

\section{Preliminaries}

We consider a spectrally negative L\'evy process $X$ which is not the negative of a subordinator or a deterministic drift with $X_{0}=0$. The passage time above a level $x\geq 0$ is defined by
\begin{equation} \label{passage}\tau_x = \inf\{t\geq 0 : X_t>x\} \;\end{equation}
which is also the right-continuous inverse of $S$ by definition. Since $S$ serves as a local time for the reflected process $S-X$ for spectrally negative L\'evy process $X$, it follows that we can express excursions of $X$ from its maximum using $\tau$. For brevity, we will use the notation $\tau$ instead of $L^{-1}$, which is the conventional notation for the right inverse of the local time.
It is well known that $\tau$ is a subordinator (killed at an exponential time if $X$ drifts to $-\infty$); see e.g. Bertoin (2007, Thm.VII.1). The points of discontinuity of $\tau=\{\tau_x:x\geq 0\}$ indicate the start of excursions $\varepsilon_x$ of $S-X$, or excursions of $X$ from its previous maximum, defined by
\[
 \varepsilon_x := \left\{ {X_{{\tau_{{x }}}}} - {X_{{\tau_{x-}} + s}}{\kern 1pt} ,{\kern 1pt} \,0 < s \le \tau_{x} - \tau_{x-}  \right\}.
 \]
for $x>0$ such that $\tau_{x} - \tau_{x-}>0$ (Kuznetsov et al., 2013).
If $\tau_{x} - \tau_{x-}=0$, take  $\varepsilon_x=\partial$ for some cemetery state $\partial$. Let ${\cal E}$ denote the space of real valued right continuous paths with left limits killed at the first hitting time of $(-\infty,0]$, endowed with the $\sigma$-algebra generated by coordinates. Then, $\{\varepsilon_x:x\geq 0, \varepsilon_x \neq \partial\}$ forms a
Poisson point process  with a characteristic measure $n$ on ${\cal E}$, which is called the excursion measure and can be constructed  as in Bertoin (2007, Chp.IV). Let $N$ denote the Poisson random measure on $(0,\infty)\times {\cal E}$ which is associated with this Poisson point process, and let $\nu$ denote its   mean measure  given by
\[
\nu(dx, d \epsilon) = dx \, n(d\epsilon)
\]
for $\epsilon \in {\cal E}$ and $x > 0$ (\c{C}inlar, 2011, Chp.6).
 Let
${\zeta}: = \inf \left\{ {s > 0:\epsilon (s) \le 0} \right\}$, which is the lifetime of an excursion $\epsilon \in {\cal E}$, and let
\[
{\bar \epsilon  }: = \sup \left\{ {{\epsilon }(s):s \le {\zeta }} \right\}\; .
\]
When $\tau_x <\infty$,   $X_{{\tau_{{x }}}}=x$ since spectrally negative L\'evy process $X$ moves continuously up having no positive jumps (Kyprianou, 2014).

For a spectrally negative L\'evy process, there exists a unique continuous increasing function $W: \mathbb{R}_+ \rightarrow \mathbb{R}_+$, called the \emph{scale function}, such that the probability that $X$ makes its first exit from an interval $[-x,y]$, $x,y>0$, at $y$ is
\begin{equation}  \label{scale}
{\p}\{I_{\tau_y} \geq -x\} = \frac{{W(x)}}{{W(x + y)}}
\end{equation}
for every $x,y>0$, and the Laplace transform of $W$ is given by
\[\int\limits_0^\infty  {{e^{ - \lambda x}}} W(x)dx = \frac{1}{{\psi (\lambda )}}   \]
see e.g. Bertoin (2007, Thm.VII.8).
The scale function $W $ is  related to excursions by
\[
\frac{W(x)}{W(y)}=\exp\left\{-\int_x^y dt\, n(\bar{\epsilon}>t)\right\}
\]
for $y>x>0$, Kyprianou (2014, Lem.8.2). It follows that
\begin{equation}  \label{nepsbar}
n(\bar{\epsilon}>t) = \frac{W'_+(t)}{W(t)}\end{equation}
where $W'_+$ denotes the right derivative of $W$. Similarly,
$n(\bar{\epsilon}\geq t) = W'_-(t)/W(t)$ where $W'_-(t)$ denotes the left derivative of $W$. The function $W$ is almost everywhere differentiable as $n(\bar{\epsilon}>t)$ has at most countably many discontinuities (Kuznetsov et al., 2013, Lem. 2.3).

\section{Joint Distribution of Maximum Loss and Maximum Gain}

We find the joint distribution of the maximum loss and the maximum gain until the time that the process passes above a level $\beta>0$, that is, till $\tau_\beta$, as defined in (\ref{passage}).  They are denoted by
\[
M_{{\tau_\beta }}^-  := \mathop {\sup }\limits_{0 \le u \le v \le {\tau_\beta }} (X_u^{} - X_v^{})\; , \quad  M_{{\tau_\beta }}^ +  = \mathop {\sup }\limits_{0 \le u \le v \le {\tau_\beta }} (X_v^{} - X_u^{}).
\]
Alternatively the following expression in terms of the excursions can be used for the maximum loss
\begin{equation}\label{kisa}
M_{{\tau_\beta }}^ -  = \mathop {\sup }\limits_{a < \beta } {\bar \varepsilon _a}\; \end{equation}
An analogous expression is valid for the maximum gain in terms of excursions of $X$ from its minimum.

The event that $X$ makes its first exit from an interval $[-x,y]$, $x,y>0$, at the upper boundary point $y$ can be written as  $\{I_{\tau_y} \geq -x\}$ and will be used to find the joint distribution of the maximum loss and the maximum gain in this section. It follows as a corollary to the following theorem, the proof of which uses similar ideas as in Salminen and Valois (2007).
\begin{thm} \label{thm2} Consider a spectrally negative L\'evy process $X$ with scale function $W$. Then, for
$u > 0,\alpha  > 0,\beta  > 0$, we have
\[
{\p}\{M_{{\tau_\beta }}^ -  \leq u,I_{\tau_{\beta}} \geq -\alpha \} = \left\{ \begin{array}{ll}
 e^{-\beta W'_+(u)/W(u)}{\rm{                         }} & u \le \alpha  \\
 \frac{{W(\alpha )}}{{W(u)}}{e^{ - (\beta  + \alpha  - u){\kern 1pt} W'_+(u)/W(u)}}{\rm{        }} & \alpha  \le u \le \alpha  + \beta  \\
 \frac{{W(\alpha )}}{{W(\alpha  + \beta )}}{\rm{  }} & \alpha  + \beta  \le u \\
 \end{array} \right.\]
In particular, $ \p\{M_{{\tau_\beta }}^ -  \leq u\} =  e^{-\beta {W'_+(u)/W(u)}} $.
\end{thm}

\prf For $u > 0,\alpha  > 0,\beta  > 0$, by (\ref{kisa}) we have
\[\{ M_{{\tau_\beta }}^ -  \leq u\}  = \{ \forall a \in (0,\beta ), {\bar \varepsilon _a} \le u, \varepsilon_a \neq \partial\}.\]
Now, consider the sets
$A=\{(a,\epsilon): 0<a<\beta, \bar{\epsilon}>u\}$ and $B=\{(a,\epsilon): 0<\alpha<\beta, \bar{\epsilon}>a+\alpha\}$.
Note that
\[{\p}\{M_{{\tau_\beta }}^ -  \leq u,I_{\tau_{\beta}} \geq -\alpha \} = {\p}\{N(A \cup B) = 0\}\] where $N$ is the Poisson random measure corresponding to the excursions.
Now, consider this probability in the following cases
\begin{enumerate}
\item if $u \leq \alpha$, then $\{ M_{{\tau_\beta }}^ -  \leq u\}$ is covered by $\{I_{\tau_{\beta}} \geq -\alpha \}$ and
\begin{eqnarray*}
\lefteqn{ \displaystyle{ {\p}\{M_{{\tau_\beta }}^ -  \leq u,I_{\tau_{\beta}} \geq -\alpha \} = {\p}\{M_{{\tau_\beta }}^ -  \leq u\} = \p\{N(A) = 0\} }}  \\
& = \exp( - \nu(A))  = {e^{ - \int\limits_A {da{\kern 1pt} {\kern 1pt} n(d\epsilon )} }} = {e^{ - \beta \,n(\bar \epsilon  > u)}} = e^{-\beta {W'_+(u)/W(u)}}
\end{eqnarray*}
by (\ref{nepsbar}), which also yields the marginal distribution.

\item if $u\ge \alpha + \beta$ then $\{ I_{\tau_{\beta}} \geq -\alpha \} $ is covered by $\{ M_{{\tau_\beta }}^ -  \leq u\}$ and we have
\[\hspace{-2cm} {\p}\{M_{{\tau_\beta }}^ -  \leq u,I_{\tau_{\beta}} \geq -\alpha \} = {\p}\{I_{\tau_{\beta}} \geq -\alpha \} = \frac{{W(\alpha )}}{{W(\alpha  + \beta )}}\]

\item if $\alpha \leq u \leq \alpha + \beta$, then we have
\begin{eqnarray*}
\lefteqn{ {\p}\{M_{{\tau_\beta }}^ -  \leq u,I_{\tau_{\beta}} \geq -\alpha \} }\\
&  = {\p}\{M_{{\tau_\beta }}^ -  \leq u,I_{\tau_{\beta}} \geq -\alpha \,|\,I_{\tau_{u-\alpha}} \geq -\alpha \}\,{\p}\{I_{\tau_{u-\alpha}} \geq -\alpha\}
\\ & = {\p_{u - \alpha }} \{ M_{{\tau_\beta }}^ -  \leq u, I_{\tau_{\beta-(u-\alpha)}} \geq -u \} \,{\p}\{I_{\tau_{u-\alpha}} \geq -\alpha\}
\end{eqnarray*}
where $\p_x$ denotes the probability law of $X$ starting from $x\in \mathbb{R}$ and we have used the strong Markov property. Similar arguments to those in i)  yield
\begin{eqnarray*}
\lefteqn{ {\p_{u - \alpha }}\{M_{{\tau_\beta }}^ -  \leq u,I_{\tau_{\beta-(u-\alpha)}} \geq -u\}}
 \\ & = {\p_{u - \alpha }}\{M_{{\tau_\beta }}^ -  \leq u\}
={e^{ - (\beta  + \alpha  - u){\kern 1pt} n(\bar \epsilon  > u)}} = e^{ - (\beta  + \alpha  - u){\kern 1pt} W'_+(u)/W(u) }
\end{eqnarray*}
where the last equality is due to (\ref{nepsbar}). Since ${\p}\{ I_{\tau_{u-\alpha}} \geq -\alpha \} =  {{W(\alpha )}}/{{W(u)}}$, the result follows.
\end{enumerate}\finish

The joint distribution of the maximum loss and the maximum gain is obtained as in the Brownian case, cf. Salminen and Vallois (2007), using the observations that $ M_{{\tau_\beta }}^ + >\beta$ and for $\alpha >0$
\[
\{M_{{\tau_\beta }}^ + - \beta \leq \alpha \} = \{ I_{\tau_{\beta}} \geq -\alpha \} \; .
\]
We put $v=\alpha +\beta$ in Theorem \ref{thm2} and get the following corollary.
\begin{cor}\label{cor1} For $v\ge \beta$ and $u\ge 0$, we have
\[\p\{M_{{\tau_\beta }}^ -  \leq u,M_{{\tau_\beta }}^ +  \leq v\} = \left\{ \begin{array}{ll}
 e^{ - \beta W'_+(u)/W(u)}{\rm{                            }}& u \le v - \beta  \\
 \frac{{W(v - \beta )}}{{W(u)}}e^{ - (v - u){\kern 1pt} W'_+(u)/W(u)}{\rm{        }}& v - \beta  \le u \le v \\
 \frac{{W(v - \beta )}}{{W(v)}}{\rm{                         }} & v \le u\, \\
 \end{array} \right.\]
\end{cor}
\vspace{5mm}

As a special case of a spectrally negative L\'{e}vy process, recall that the scale function for Brownian motion with drift $\mu$ and $\sigma^2=1$ is given by $$S^{\mu}(x):=W(x)=\frac{1}{2\mu}(1-e^{-2\mu x}) \; .  $$
Then, we get
$$W^{'}_{+}(u)/W(u)=2\mu/(e^{2\mu u}-1)=1/S^{-\mu}(u).$$ By substituting this term and the scale function $S^\mu$ in the statements of Theorem \ref{thm2} and in Corollary \ref{cor1}, we see that the results for the Brownian motion with drift  given in Salminen and Vallois (2007, Prop. 2.1, Cor. 2.2) are recovered.

\section{Marginal Distributions up to an Exponential Time}

Let the right inverse of $\psi$ be denoted with $\phi$, which is given by
\[
\phi(\gamma)=\sup\{\lambda\geq 0: \psi (\lambda)=\gamma\}
\]
for $\gamma>0$. In addition to the scale function $W$, there exist a family of increasing functions $W^{(q)}: \mathbb{R}_+ \rightarrow \mathbb{R}_+$ which satisfy
\begin{equation}  \label{laplace}
\int\limits_0^\infty  {{e^{ - \lambda x}}} {W^{(q)}}(x)dx = \frac{1}{{\psi (\lambda ) - q}}{\rm{ ,      }}
\end{equation}
for $q \geq 0$ and $\lambda > \phi(q)$, and also functions $Z^{(q)}$ defined on $\mathbb{R}$ by
\begin{equation}  \label{Z}
{Z^{(q)}}(x) = 1 + q\int_0^x {{W^{(q)}}(y)dy} \; .
\end{equation}
The properties of the so-called $q$-scale functions $W^{(q)}$,  and $Z^{(q)}$ can be gathered from Kyprianou (2014, Thm.8.1). Note that $ W^{(0)}\equiv W$ above. Let $W_ + ^{(q)'}(x)$ denote the right derivative of ${W^{(x)}}$ at $a$.

 \begin{prop}\label{prop1} For the maximum loss and the maximum gain up to an independent exponential time $T$ with parameter $\gamma>0$
 \begin{eqnarray*}
\p\{M_T^- >a\}& = & {Z^{(\gamma)}}(a) - \gamma\, {[{W^{(\gamma)}}(a)]^2}/W_ + ^{(\gamma)'}(a) \\
\p\{M_T^+ >a\}& = & \frac{{Z^{(\gamma)}}(0)}{{Z^{(\gamma)}}(a)}
 \end{eqnarray*}
for $a>0$.
 \end{prop}

 \prf For the loss process, the first time to pass the level $a$ is defined by
\[{\overline{\sigma}_a} = \inf \{ v \ge 0:{S_v} - {X_v} > a\} \]
 The Laplace transform of $\overline{\sigma}_a$ has been found in Avram et al. (2004, Thm.1) as
\[ \
\E({e^{ - \gamma{\overline{\sigma} _a}}}) = {Z^{(\gamma)}}(a) - \gamma{[{W^{(\gamma)}}(a)]^2}/W_ + ^{(\gamma)'}(a) \;.
\]
Clearly, we have $\p\{M_T^- >a\}=\E(e^{-\gamma \overline{\sigma}_a})$. Similarly, we can use the  formula for  $\E({e^{ - \gamma{\underline{\sigma}\, _a}}})$ from Pistorius (2004, Prop.2) to get the distribution of $M_T^+$, where $\underline{\sigma}\, _a=\inf\{v\geq 0 : X_v-I_v >a\}$. \finish

\vspace{5mm}

In Hubalek and Kyprianou (2010), the scale function $W^{(\gamma)},$ $\gamma\geq0$, for Brownian motion with drift $\mu$ and $\sigma^2=1$ is given for $x\geq 0$ as
$$W^{(\gamma)}(x) = \frac{2}{\sqrt{2\gamma+\mu^{2}}}e^{-\mu x}\sinh(x\sqrt{2\gamma+\mu^{2}})$$ which can be obtained by inverting the Laplace transform, Equation (\ref{laplace}). Using this $\gamma$-scale function  and Equation (\ref{Z}), we have $Z^{(\gamma)}(0)=1$, and we obtain
\begin{eqnarray*}Z^{(\gamma)}(a)&=& 1 + \gamma\int_0^a {\frac{2}{\sqrt{2\gamma+\mu^{2}}}e^{-\mu x}\sinh(x\sqrt{2\gamma+\mu^{2}})dx} \\
&=&\frac{\mu}{\sqrt{2\gamma+\mu^{2}}}e^{-a\mu}\sinh(a\sqrt{2\gamma+\mu^{2}})+e^{-a\mu}\cosh(a\sqrt{2\gamma+\mu^{2}})
 \end{eqnarray*}
and
$$\frac{W^{(\gamma)}(a)^2}{W_ + ^{(\gamma)'}(a)}=\frac{2e^{-\mu a}\sinh^{2}(a\sqrt{2\gamma+\mu^{2}})}{(2\gamma+\mu^{2})\cosh(a\sqrt{2\gamma+\mu^{2}})-\mu\sqrt{2\gamma+\mu^{2}}\sinh(a\sqrt{2\gamma+\mu^{2}})} \; .$$
By substituting these expressions in Proposition \ref{prop1}, the marginal distributions for Brownian motion with drift case are found as
 \begin{eqnarray*}
\p\{M_T^- >a\}& = & \frac{1}{e^{a\mu}\cosh(a\sqrt{2\gamma+\mu^{2}})-\frac{\mu}{\sqrt{2\gamma+\mu^{2}}}e^{a\mu}\sinh(a\sqrt{2\gamma+\mu^{2}})} \\
\p\{M_T^+ >a\}& = & \frac{1}{e^{-a\mu}\cosh(a\sqrt{2\gamma+\mu^{2}})+\frac{\mu}{\sqrt{2\gamma+\mu^{2}}}e^{-a\mu}\sinh(a\sqrt{2\gamma+\mu^{2}})}
 \end{eqnarray*}
for $a>0$. Note that, these results match with the formulas given in Eqns. (1.8) and (1.9) in Salminen and Vallois (2007).

\vspace{0.5cm}
\noindent{\bf Acknowledgement.}
The authors would like to express their gratitude to an anonymous reviewer for detailed comments and guidance on the presentation of their results.

\vspace{0.5cm}


\begin{thebibliography}{}

\bibitem{Avra} F. Avram,  A.E. Kyprianou, M.R. Pistorius: Exit Problems for Spectrally Negative L\'evy Processes and Applications to (Canadized) Russian Options, Annals of Applied Probability, 14: 215-238 (2004).

\bibitem{bert} J. Bertoin: L\'evy Processes, Cambridge University Press, Cambridge (2007).

\bibitem{cinl} E. \c{C}inlar: Probability and Stochastics, Springer, New York (2011).

\bibitem{Hubalek} F. Hubalek,  A.E. Kyprianou: Old and New Examples of Scale Functions
for Spectrally Negative L\'evy Processes, Progress in Probability, Vol. 63, 119-145 (2010).

\bibitem{KuztKypr}A. Kuznetsov, A.E. Kyprianou, V. Rivero: The Theory of Scale Functions for Spectrally Negative L\'evy Processes,
L\'evy Matters II, Lecture Notes in Mathematics, 2061: 97-186 (2013).

\bibitem{Kypr} A.E. Kyprianou: Fluctuations of L\'evy Processes with Applications, 2nd Ed., Springer, New York (2014).

\bibitem{MijoPist} A. Mijatovic, M.R. Pistorius: On the Drawdown of Completely Asymmetric L\'evy Processes, Stoch. Proc. Appl.,  122: 3812-3836 (2012).

\bibitem{pist} M.R. Pistorius: On Exit and Ergodicity of the Spectrally One-Sided
L\'evy Process Reflected at Its Infimum, Journal of Theoretical Probability, 17: 183-220 (2004).

\bibitem{salm} P. Salminen, P. Vallois: On maximum increase and decrease of Brownian motion, Ann. I. H. Poincar\'e, 43: 655-676 (2007).

\bibitem{ceren}
C. Vardar-Acar, C. L. Zirbel, G.J. Sz\'ekely: On the correlation of the supremum and the infimum and of maximum gain and maximum loss of Brownian motion with drift, Journal of Computational and Applied Mathematics, 248: 61-75 (2013).

\end{thebibliography}
\end{document}